\input amstex.tex
\documentstyle{amsppt}
\magnification=1200

\loadbold

\def\Co#1{{\Cal O}_{#1}}
\def\Fi#1{\Phi_{|#1|}}
\def\fei#1{\phi_{#1}}

\def\lrw{\longrightarrow}

\def\Bbbp1{{\Bbb P}^1}
\def\simlin{\sim_{\text{\rm lin}}}
\def\simnum{\sim_{\text{\rm num}}}
\def\Div{\text{\rm Div}}

\def\dim{\text{\rm dim}}
\def\ulrcorner#1{\ulcorner{#1}\urcorner}
\def\llrcorner#1{\llcorner{#1}\lrcorner}
\def\Bbbq{\Bbb Q}
\def\kod{\text{\rm kod}}

\TagsOnRight

\topmatter
\title
On canonically derived families of surfaces of general type over curves
\endtitle
\author Meng Chen
\endauthor
\thanks *Partially supported by the {\bf National Natural Science Foundation
of China}
\endthanks
\endtopmatter

\document
\leftheadtext{M. Chen}
\rightheadtext{Canonically derived families of surfaces of general type}

\centerline {\it Department of Applied Mathematics, Tongji University}
\centerline {\it Shanghai, 200092, P. R. China} 
\centerline {{\it E-mail}: mchen\@mail.tongji.edu.cn}
\bigskip

The aim of this paper is to study a family of surfaces of general type over a curve which is derived from a multicanonical system of a smooth projective threefold of general type. As far as we know, there are no systematical references yet on this topic. In fact, to study this kind of families is an important step of the classification theory.

We always suppose that the ground field is algebraically closed of characteristic zero. Let $X$ be a smooth projective variety of general type with dimension $d$. We say that $|mK_X|$ {\it is composed of a pencil of varieties of dimension $d-1$} if $\dim_{\Bbb C}H^0(X,\Co{X}(mK_X))\ge 2$ and $\dim\fei{m}(X)=1$, where $m$ is a positive integer and $\fei{m}:=\Fi{mK_X}$ is the rational map defined by the system $|mK_X|$. Set $P_m(X):=\dim_{\Bbb C} H^0(X, \Co{X}(mK_X))$. We call $P_m(X)$ {\it the m-th genus of $X$}, which is an important birational invariant.  

Now suppose $|mK_X|$ is composed of a pencil.
Take possible blow-ups $\pi:X'\lrw X$, according to Hironaka, such that $g_m:=\fei{m}\circ\pi$ is a morphism onto its image. Denote 
$$W_m:=\overline{\fei{m}(X)}\subset {\Bbb P}^{P_m(X)-1}$$
and let 
$$g_m: X'\overset{f_m}\to\lrw C\overset{\psi_m}\to\lrw W_m$$
be the Stein factorization of $g_m$, where $C$ is a smooth projective curve of genus $b:=g(C)$. Then we have the following commutative diagram: 
$$\CD
X' @>{id.}>> X' @>\pi>>  X\\
@V{f_m}VV   @VV{g_m}V    @VV\fei{m}V\\
C  @>>\psi_m> W_m  @<<{id.}< W_m
\endCD$$
where we note that $\fei{m}$ is only a rational map.
Denote by $F$ a general fiber of $f_m$. Then $F$ is a smooth projective variety of general type of dimension $d-1$.
We say that the fibration $f_m:X'\lrw C$ is {\it a derived family from the m-canonical pencil $|mK_X|$}, which is the main object of this paper. 

When $d=2$ and $m=1$, there are infinite number of such families according to \cite{Be} and the classification is still incomplete according to \cite{Ca2}
where Catanese first constructed a canonically derived family of curves with non-constant moduli. When $d=2$ and $m\ge 2$, only one possibility occurs according to Theorem 1 of \cite{X1}. Explicitly,  such a family of curves is derived from the bicanonical system of a smooth surface $S$ of general type with the invariants $(K^2, p_g)=(1,0)$, where $K^2$ and $p_g$ are both the invariants of the minimal model of $S$.  We note that, in this situation, $P_2(S)=2$ and both $p_g(S)$ and $q(S)$ take minimal values.

It is natural we turn our interest to higher dimensional case. Here, we only treat the case $d=3$. The existence of canonically derived family of surfaces is undoubted. One might refer to \cite{R} for many examples with $P_m=2$. However it isn't quite clear about the bulk of the set of these families. 
This paper aims to build some basic facts and to study these families in terms of birational invariants of  the total space as well as those of a general fiber. We observed that Koll\'ar proved the following result.

\proclaim{Theorem 0} (Theorem 6.1 of \cite{Ko1}) Let $X$ be a smooth projective 3-fold of general type. If $m\ge 3$ and $\dim\fei{m}(X)\le 2$, then $q(X):=\dim H^1(X,\Co{X})\le 3$.
\endproclaim

As was pointed out by Koll\'ar, one can get $q(X)=0$ whenever $m\gg 0$ under
the condition of Theorem 0. Unfortunately, the bound of $m$ would be much
bigger by virtue of his method. Since, in our case, $\dim\fei{m}(X)=1$, we
should be able to get  more explicit information even if $m$ is small. On
the basis of a detailed classification (Theorems 3.3, 3.4), we obatained the
following results in this paper.

\proclaim{Theorem 1} Let $f: X\lrw C$ be a derived family of surfaces from the m-canonical pencil $|mK_X|$ of a smooth projective 3-fold $X$ of general type. 
Suppose $m\ge 3$. 
Then $C$ is either an elliptic curve or ${\Bbb P}^1$ and $f$
has the following properties.

(i) Either $q(X)\le 1$ or $p_g(X)\le 1$. 

(ii) $q(X)\le 2$ whenever $m\ge 11$.

(iii) $q(X)\le 2$ whenever $m\ge 7$ and $p_g(X)>0$.

(iv) $p_g(X)\le 1$ whenever $m\ge 7$. $p_g(X)\le 2$ whenever $5\le m\le 6$. $p_g(X)\le 3$ whenever $m=4$. $p_g(X)\le 5$ whenever $m=3$.

(v) If $q(X)=3$, then either $P_m(X)=2$  or $\dim\fei{m+1}(X)\ge 2.$
\endproclaim

\proclaim{Theorem 2} Let $f: X\lrw C$ be a derived family of surfaces from the bicanonical pencil $|2K_X|$ of a smooth projective 3-fold $X$ of general type. 
Then $C$ is either an elliptic curve or ${\Bbb P}^1$ and $f$
has the following properties.

(i) Either $q(X)\le 2$ or $p_g(X)\le 2$.

(ii) If $q(X)\ge 3$, then either $P_2(X)=2$ or $\dim\fei{3}(X)\ge 2$.
\endproclaim

In the final section, we would like to give an appendix to Koll\'ar's method on how to determine the bounds of $m$ so as to get $q(X)\le 1$. The bounds obtained are much better than that of Koll\'ar. However, we feel that they are still far from being the optimal ones. The result is as follows.   

\proclaim{Corollary 3} Let $f: X\lrw C$ be a derived family of surfaces from the m-canonical pencil $|mK_X|$ of a smooth projective 3-fold $X$ of general type. 
Then  

(i)  $q(X)\le 1$ whenever $m\ge 82$.

(ii) $q(X)=0$ whenever $m\ge 143$.
\endproclaim

\head 1. Preliminaries \endhead

\subhead 1.1 Convention\endsubhead
Let $X$ be a normal projective variety of dimension $d$. We denote by
$\Div(X)$
the group of Weil divisors on $X$. An element $D\in\Div(X)\otimes{\Bbb Q}$
is called a ${\Bbb Q}$-{\it divisor}. A ${\Bbb Q}$-divisor $D$ is
said to be ${\Bbb Q}$-{\it Cartier} if $mD$ is a Cartier divisor for some
positive integer $m$. For a ${\Bbb Q}$-Cartier divisor
$D$ and an irreducible curve $C\subset X$, we can define the intersection
number $D\cdot C$ in a natural way. A $\Bbbq$-Cartier divisor $D$ is called
{\it nef} (namely {\it numerically effective}) if $D\cdot C\ge 0$ for any
effective curve
$C \subset X$. A nef divisor $D$ is called {\it big} if $D^d>0$. We say that
$X$ is $\Bbbq$-{\it factorial} if every Weil divisor on $X$ is
$\Bbbq$-Cartier.
For a Weil divisor $D$ on $X$, write $\Co{X}(D)$ as the corresponding
reflexive sheaf.
Denote by $K_X$ a canonical divisor of $X$, which is a Weil divisor. $X$ is
called
{\it minimal} if $K_X$ is a nef $\Bbbq$-Cartier divisor. $X$ is said to be
of general type if $\kod(X)=\dim(X)$. For a positive integer $m$, we set
$\omega_X^{[m]}:=\Co{X}(mK_X).$
We use \cite{R} as a nice reference for the definition of {\it canonical, terminal singularities}. According to both \cite{KMM} and \cite{K-M}, any given smooth projective 3-fold $Y$ of general type has a minimal model $X$ which has only ${\Bbb Q}$-factorial terminal singularities.   

\subhead 1.2 Vanishing theorems\endsubhead
Let $D=\sum a_iD_i$ be a ${\Bbb Q}$-divisor on $X$ where the $D_i's$ are distinct prime divisors and $a_i\in {\Bbb Q}$. We define

{\it the round-down} $\llrcorner{D}:=\sum \llrcorner{a_i}D_i$, where $\llrcorner{a_i}$ is the integral part of $a_i$.

{\it the round-up} $\ulrcorner{D}:=-\llrcorner{-D}$.

{\it the fractional part} $\{D\}:=\ulrcorner{D-\llrcorner{D}}$.

Throughout this paper, we will use the Kawamata-Viehweg vanishing theorem (\cite{Ka1}, \cite{KMM} and \cite{V}) in the following forms.
\proclaim{Theorem 1.1}
Let $X$ be a smooth complete variety, $D\in \Div(X)\otimes{\Bbb Q}$. Assume the following two conditions: 

(i) $D$ is nef and big;

(ii) the fractional part of $D$ has supports with only normal crossings.

\noindent Then $H^i(X,\Co{X}(K_X+\ulrcorner{D}))=0$ for all $i>0$.
\endproclaim

\proclaim{Theorem 1.2} Let $X$ be a normal projective variety with only canonical singularities. Let $D$ be a ${\Bbb Q}$-Cartier Weil divisor such that $D$ is nef and big. Then $H^i(X,\Co{X}(K_X+D))=0$ for all $i>0$.
\endproclaim

\subhead 1.3 Semi-positivity \endsubhead Let $C$ be a smooth projective curve and ${\Cal E}$ be a vector bundle on $C$. We call 
$$\mu({\Cal E}):=\frac{\deg({\Cal E})}{\text{rk}({\Cal E})}$$
{\it the slope of} ${\Cal E}$.
According to \cite{H-N}, there is the Harder-Narasimhan filtration
$$0={\Cal E}_0\subset {\Cal E}_1\subset\cdots\subset{\Cal E}_{n-1}\subset{\Cal E}_n={\Cal E},$$
where the quotient ${\Cal E}_i/{\Cal E}_{i-1}$ is a semistable vector bundle and $$ \mu({\Cal E}_i/{\Cal E}_{i-1})>\mu({\Cal E}_{i+1}/{\Cal E}_i)$$
for all $i$. We define
$$\mu_{min}({\Cal E}):=\mu({\Cal E}/{\Cal E}_{n-1}),$$
which is called {\it the minimal slope of} ${\Cal E}$.

\definition{Definition 1.3} The vector bundle ${\Cal E}$ is said to be {\it semi-positive} if $\mu_{min}({\Cal E})\ge 0$.
\enddefinition
According to \cite{Ka2}, \cite{Ko2}, \cite{N} and \cite{O}, we have the following

\proclaim{Fact 1.4} Let $X$ be a smooth projective 3-fold and $f:X\lrw C$ be a proper morphism with connected fibers onto a smooth projective curve $C$. Then
both $f_*\omega_{X/C}^{\otimes m}$ and $R^if_*\omega_{X/C}$ are semi-positive vector bundles on $C$ for all $m>0$ and $i>0$. In particular, 
$$\deg f_*\omega_{X/C}^{\otimes m}\ge 0\ \text{and}\ 
\deg R^if_*\omega_{X/C}\ge 0.$$
\endproclaim

\subhead 1.4 Basic formulae \endsubhead
Let $X$ be a smooth projective 3-fold and $f:X\lrw C$ be a fibration onto the smooth projective curve $C$. Denote $b:=g(C)$. From the spectral sequence
$$E_2^{p,q}:=H^p(C, R^qf_*\omega_X)\Longrightarrow E^n:=H^n(X,\omega_X),$$
one obtains the following formulae
$$q(X):=h^1(X,\Co{X})=b+h^1(C,R^1f_*\omega_X)\tag{1.1}$$
$$h^2(\Co{X})=h^1(C, f_*\omega_X)+h^0(C,R^1f_*\omega_X).\tag{1.2}$$

\head 2. Lemmas \endhead

\proclaim{Lemma 2.1} Let $X$ be a smooth projective 3-fold of general type. $m\ge 2$ is an integer. Suppose that $|mK_X|$ is composed of a pencil of surfaces. Keep the same notations as in the first page of this paper. We have a derived fibration $f_m:X'\lrw C$. Then $C$ is either an elliptic curve or 
${\Bbb P}^1$.
\endproclaim
\demo{Proof}
Suppose $b>0$. Then $\fei{m}$ is a morphism. We have a derived fibration
$$f:=f_m: X\lrw C.$$
Let ${\Cal E}_0$ be the saturated sub-bundle of $f_*\omega_X^{\otimes m}$ which is generated by $H^0(C, f_*\omega_X^{\otimes m})$. Since $|mK_X|$ is composed of a pencil and $\fei{m}$ factors through $f$, ${\Cal E}_0$ should be a line bundle on $C$. Denote ${\Cal E}:=f_*\omega_X^{\otimes m}$. Then we have the following extension
$$0\lrw {\Cal E}_0\lrw {\Cal E}\lrw {\Cal E}_1\lrw 0$$
and the exact sequence
$$f_*\omega_{X/C}^{\otimes m}\lrw {\Cal E}_1\otimes \omega_C^{\otimes -m}\lrw 0.$$
Note that $r:=\text{rk}({\Cal E})=h^0(F, mK_F)\ge 2$ because the general fiber $F$ is a smooth projective surface of general type. According to Fact 1.4, $f_*\omega_{X/C}^{\otimes m}$ is semi-positive. Therefore $\deg({\Cal E}_1\otimes\omega_C^{\otimes -m})\ge 0$, i.e.
$$\deg{\Cal E}_1\ge 2m(r-1)(b-1).$$
We have
$$\align
h^1({\Cal E}_0)&\ge h^0({\Cal E}_1)\ge \deg{\Cal E}_1+(r-1)(1-b)\\
&\ge (2m-1)(r-1)(b-1).   \tag{2.1}
\endalign$$
Suppose $h^1({\Cal E}_0)>0$. Since $\deg {\Cal E}_0>0$, according to the Clifford's theorem, we have
$$\deg {\Cal E}_0\ge 2h^0({\Cal E}_0)-2\ge h^0({\Cal E}_0)=P_m(X)\ge 2.$$
On the other hand, we have 
$$h^1({\Cal E}_0)=h^0({\Cal E}_0)-\deg({\Cal E}_0)+b-1\le b-1.\tag{2.2}$$
Thus, by (2.1) and (2.2), we get 
$$b-1\ge (2m-1)(r-1)(b-1).$$
The only possibility is $b=1$. When $h^1({\Cal E}_0)=0$, we also automatically have $b=1$ from (2.1). The proof is complete.
\qed
\enddemo

\proclaim{Lemma 2.2} Let ${\Cal E}$ be a vector bundle of rank $r$ on a smooth projective curve $C$. Suppose ${\Cal E}\otimes \omega_C^{-1}$ is semi-positive.
Then we have $h^1(C,{\Cal E})\le r.$
\endproclaim
\demo{Proof}
Suppose there are $r+1$ independant sections
$$s_1, s_2, \cdots, s_{r+1}\in H^0(C, ({\Cal E}\otimes\omega_C^{-1})^\lor)\cong
H^1(C,{\Cal E}).$$
Denote ${\Cal E}':=({\Cal E}\otimes\omega_C^{-1})^\lor$.
For any point $x\in C$, the stalk ${\Cal E}_x'$ is an $\Co{C,x}$-module of rank $r$. This means that $s_{1,x}$, $s_{2,x}$, $\cdots$, $s_{r+1,x}$ are algebraically dependant in ${\Cal E}'_x$. Thus there are $r+1$ nontrivial germs
$$f_1, f_2, \cdots, f_{r+1}\in \Co{C,x}$$
such that $\sum_{i=1}^{r+1}f_i(x)s_{i,x}(x)=0$. Now set 
$$s:=\sum_{i=1}^{r+1} f_i(x)s_i\in H^0(C,({\Cal E}\otimes\omega_C^{-1})^\lor).$$
$s$ is a non-zero section. Otherwise $s_1$, $\cdots$, $s_{r+1}$ are dependant.
Because $s$ vanishes at $x$, $s$ defines a line bundle ${\Cal L}$ which has positive degree. So ${\Cal E}\otimes\omega_C^{-1}$ has a quotient bundle with negative degree. This contradicts to the semi-positivity of ${\Cal E}\otimes\omega_C^{-1}$. The proof is completed.
\qed\enddemo

\proclaim{Corollary 2.3} Let $f: X\lrw C$ be a fibration from a smooth projective 3-fold $X$ onto a smooth projective curve $C$. Let $F$ be a general fibre of $f$ and set $b:=g(C)$. Then $q(X)\le b+q(F)$.
\endproclaim
\demo{Proof} This is a direct result from Fact 1.4, (1.1) and Lemma 2.2.
\enddemo

\head 3. Proof of the main theorems \endhead

Since the behavior of pluricanonical maps is birationally invariant, we may suppose that $X$ is a normal projective minimal 3-fold of general type with only ${\Bbb Q}$-factorial terminal singularities. We make this assumption so as to utilize vanishing theorems. Now suppose that $|mK_X|$ is composed of a pencil of surfaces. We can use the same set up as in the first page of this paper. An extra point is that we can take the modification $\pi_m:X'\lrw X$ such that 
$\pi_m^*(mK_X)$ has supports with only normal crossings. 
We keep the same notations. Then we get a derived fibration $f_m:X'\lrw C$. Denote by $F$ a general fiber of $f_m$ and by $b$ the genus of $C$. Sometimes we simply denote $f_m$ by $f$ and $\pi_m$ by $\pi$, respectively. 

\proclaim{Proposition 3.1}  Let $X$ be a normal projective minimal 3-fold of general type with only ${\Bbb Q}$-factorial terminal singularities. Suppose $|mK_X|$ is composed of a pencil of surfaces. If $p_g(X)>0$, then $p_g(F)=1$ under one of the following conditions.

(1) $m=2$, $b=0$ and $P_m(X)\ge 4$.

(2) $m\ge 3$, $b=0$ and $P_m(X)\ge 3$.

(3) $m=2$, $b=1$ and $P_m(X)\ge 3$.

(4) $m\ge 3$, $b=1$ and $P_m(X)\ge 2$.
\endproclaim
\demo{Proof}
Because $p_g(X')=p_g(X)>0$, we can choose an effective divisor $D_1\in |K_{X'}|$.
We write
$$K_{X'}=\pi^*(K_X)+\sum a_iE_i,$$
where $a_i\in {\Bbb Q}^{+}$, $E_i$ is an exceptional prime divisor for all $i$. We note that
$$\pi^*(K_X)=K_{X'}-\sum a_iE_i=D_1-\sum a_iE_i$$
is an effective ${\Bbb Q}$-divisor. So we have 
$$\pi^*(K_X)=\sum b_iG_i=G_v+G_h,$$ 
where $b_i\in {\Bbb Q}^{+}$, $G_i$ is a prime divisor on $X'$ for all $i$, the support of $G_v$ is contained in fibers of $f$ and $G_h$ is the horizontal part of $\pi^*(K_X)$. Both $G_v$ and $G_h$ are effective ${\Bbb Q}$-divisors. So we have 
$$(m-1)\pi^*(K_X)=(m-1)G_v+(m-1)G_h$$
$$m\pi^*(K_X)=mG_v+mG_h.$$
Suppose $M_m$ is the movable part of $|mK_{X'}|$. Then $M_m\le_{\Bbb Q} mG_v$ because the support of $M_m$ is vertical. Thus we can write
$$mG_v=M_m+G_v'=M_m+\sum c_iG_i',$$ 
where $c_i\in {\Bbb Q}^{+}$, $G_i'$ is a prime divisor on $X'$ for all $i$. Therefore 
$$m\pi^*(K_X)=M_m+(G_v'+mG_h),$$
where $G_v'+mG_h$ is an effective ${\Bbb Q}$-divisor. We can suppose
$$M_m\simlin\sum_{i=1}^a F_i\simnum aF,$$
where 
$$a=\cases P_m(X)-1, & \text{if}\ b=0,\\
P_m(X), &\text{if}\ b=1.\endcases $$
Thus we have
$$\pi^*(K_X)\simnum\frac{a}{m}F+\frac{1}{m}G_v'+G_h$$
$$(m-1)\pi^*(K_X)\simnum\frac{a(m-1)}{m}F+\frac{m-1}{m}G_v'+(m-1)G_h.$$
Denote $a':=\frac{a(m-1)}{m}$. We can see that $a'>1$ under one of the assumptions within (1) through (4) of the proposition. So
$$(m-1)\pi^*(K_X)-F-\frac{m-1}{ma'}G_v'-\frac{m-1}{a'}G_h\simnum
(m-1)(1-\frac{1}{a'})\pi^*(K_X)$$
is nef and big and its fractional part has supports with only normal crossings.
According to the Kawamata-Viehweg vanishing theorem, we get 
$$H^1(X',K_{X'}+\ulrcorner{(m-1)\pi^*(K_X)-\frac{m-1}{ma'}G_v'-\frac{m-1}{a'}G_h}-F)=0.$$
Set $G'':=\ulrcorner{(m-1)\pi^*(K_X)-\frac{m-1}{ma'}G_v'-\frac{m-1}{a'}G_h}$.
Then $G''\le\ulrcorner{(m-1)\pi^*(K_X)}$ and so
$$K_{X'}+G''\le K_{X'}+\ulrcorner{(m-1)\pi^*(K_X)}.$$
Therefore we have 
$$\dim\Fi{K_{X'}+G''}(F)=0$$ 
for a general fiber $F$. {}From the exact sequence
$$0\lrw\Co{X'}(K_{X'}+G''-F)\lrw\Co{X'}(K_{X'}+G'')\lrw\Co{F}(K_F+G''|_F)\lrw 0,$$
we get the surjective map
$$H^0(X', K_{X'}+G'')\lrw H^0(F, K_F+G''|_F).$$
This means
$$|K_{X'}+G''|\bigm|_F=\bigm|K_F+G''|_F\bigm|.$$
Noting that
$$G''|_F=\ulrcorner{(m-1)\pi^*(K_X)-\frac{m-1}{ma'}G_v'-\frac{m-1}{a'}G_h}|_F
=\ulrcorner{(m-1)(1-\frac{1}{a'})G_h}|_F$$
is an effective divisor, we have $h^0(F, K_F+G''|_F)\ge 2$ whenever $p_g(F)\ge 2$. This would lead to $\dim\fei{m}(F)\ge 1$, which is impossible. So we should have $p_g(F)=1$ because $p_g(F)>0$ under the assumption $p_g(X)>0$.
\qed\enddemo

\remark{Remark 3.2} The assumption $p_g(X)>0$ in Proposition 3.1 is important. If $p_g(X)=0$, the above method is invalid because we don't know whether $G''|_F$ is effective.
\endremark

\proclaim{Theorem 3.3} Let $f: X\lrw C$ be a derived family of surfaces from the m-canonical pencil $|mK_X|$ of a smooth projective 3-fold $X$ of general type. 
Let $F$ be a general fiber of $f$ and denote $b:=g(C)$.
Suppose $m\ge 3$. Then one of the following occurs:

(A1) $m=5,\ 6$, $p_g(X)=2$, $q(X)=0$, $b=0$ and $p_g(F)=q(F)=1$.

(A2) $m=4$, $2\le p_g(X)\le 3$, $q(X)=0$, $b=0$ and $p_g(F)=1$.
 
(A3) $m=3$, $2\le p_g(X)\le 5$, $q(X)\le 1$, $b=0$ and $p_g(F)=1$.

(B0) $p_g(X)=1$, $q(X)\le b+1$ and $p_g(F)=1$.

(B1) $p_g(X)=P_2(X)=\cdots=P_{m-1}(X)=1$, $P_m(X)=2$, $\dim\fei{m+1}(X)\ge 2$, $b=0$ and $p_g(F)\ge 2$.

(B2) $p_g(X)=P_2(X)=\cdots=P_{m-1}(X)=1$, $P_m(X)=P_{m+1}(X)=2$, $\dim\fei{m+2}(X)\ge 2$, $b=0$ and $p_g(F)\ge 2$.

(B3) $p_g(X)=P_2(X)=\cdots=P_{m-2}(X)=1$, $P_{m-1}(X)=P_m(X)=2$, $\dim\fei{m+1}(X)\ge 2$, $b=0$ and $p_g(F)\ge 2$.

(C0) $p_g(X)=0$, $q(X)\le 2$ and $q(F)\le 2$.

(C1) $p_g(X)=0$, $\dim\fei{m+1}(X)\ge 2$ and $q(F)\ge 2$.

(C2) $p_g(X)=0$, $P_m(X)=2$, $\dim\fei{m+2}(X)\ge 2$, $b=1$ and $q(F)=2$.

(C3) $p_g(X)=0$, $P_m(X)=2$, $\dim\fei{2m}(X)\ge 2$, $b=0$ and $q(F)\ge 3$.

(C4) $p_g(X)=0$, $P_{m}(X)=2$, $P_{2m}(X)=3$, $\dim\fei{2m+1}(X)=3$, $b=0$ and $q(F)\ge 3$.
\endproclaim
\demo{Proof} We formulate the proof through three steps. Though the proof is slightly longer, it's a case by case discussion.

Step 1. $p_g(X)\ge 2$ 

Suppose $b=1$. In this situation, we see that the movable part of $|3K_X|$ defines a morphism. Because $p_g(X)\ge 2$, $\dim\fei{1}(X)=1$ and both $\fei{1} $ and $\fei{3}$ derive the same fibration $f:X\lrw C$. So the movable part of $|K_X|$ is also base point free. Let $M_1$ be the movable part of $|K_X|$. Then $M_1\simlin \sum F_i$. Let $F$ be a general fiber of $f$. Because the singularities on $X$ are all isolated, $F$ is a smooth projective surface of general type. By Theorem 2.2, we see that $H^1(X, 2K_X)=0$. Therefore we have 
$$|2K_X+\sum F_i||_F=|2K_F|.$$
Because $p_g(F)>0$, $\Fi{2K_F}$ is generically finite by Theorem 1 of \cite{X1}. This means that $\fei{3}$ is generically finite and so is $\fei{m}$ for $m\ge 4$. So we only have to consider the case when $b=0$.

Now we have a fibration $f:X'\lrw {\Bbb P}^1$. Because $p_g(X)\ge 2$, we have $P_3(X)\ge 4$. Thus, by Proposition 3.1, we see that $p_g(F)=1$.  In this situation, Koll\'ar's technique (the proof of Corollary 4.8 in \cite{Ko1}) is still effective. Let $p_g(X)=k+1$, $k\ge 1$. Since $|K_{X'}|$ is composed of a pencil, we have ${\Cal O}(k)\hookrightarrow f_*\omega_{X'}$ on ${\Bbb P}^1$. If $k\ge 5$, then 
$${\Cal O}(5)\hookrightarrow {\Cal O}(k)\hookrightarrow f_*\omega_{X'}.$$
Thus we have
$${\Cal E}:={\Cal O}(1)\otimes f_*\omega_{X'/{\Bbb P}^1}^2={\Cal O}(5)\otimes f_*\omega_{X'}^2\hookrightarrow f_*\omega_{X'}^3.$$
The local sections of $f_*\omega_{X'}^2$ give the bicanonical map of the fiber $F$ and they extend to global sections of ${\Cal E}$, because ${\Cal E}$ is generated by global sections. On the other hand, $H^0({\Bbb P}^1, {\Cal E})$ can distinguish different fibers of $f$ because $f_*\omega_{X'/{\Cal P}^1}^2$ is a sum of line bundles with nonnegative degree on ${\Bbb P}^1$. So $H^0({\Bbb P}^1, {\Cal E})$ gives a generically finite map on $X'$ and so does $H^0(X', 3K_{X'})$.
This contradicts to our assumption of $\dim\fei{3}(X)=1$. Thus we have $k\le 4$, i.e. $p_g(X)\le 5$. By virtue of this technique, we have 
$$\align
k=3,\ 4,\  &{\Cal E}\hookrightarrow f_*\omega_{X'}^4\Longrightarrow m=3\Longrightarrow\  (A3)\\
k=2, \ \ \ \  \ &{\Cal E}\hookrightarrow f_*\omega_{X'}^5\Longrightarrow m\le 4\Longrightarrow\  (A2),\ (A3)\\
k=1, \ \ \  \ \ &{\Cal E}\hookrightarrow f_*\omega_{X'}^7\Longrightarrow m\le 6\Longrightarrow\  (A1),\  (A2),\ (A3).
\endalign$$
In order to complete the proof for this case, we have to prove $q(X)=0$ for (A1), (A2) and that the only possibility of $F$ in (A1) is $p_g(F)=q(F)=1$. 
Suppose $q(X)=1$ in cases (A1) and (A2). Then $q(F)=1$ and $R^1f_*\omega_{X'}\cong \omega_{{\Bbb P}^1}$. Because $f_*\omega_{X'}$ is of positive degree, we have 
$h^1({\Bbb P}^1, f_*\omega_{X'})=0$. So by (1.2), $h^2(\Co{X'})=0$. Then we have $\chi(\Co{X'})\le -2$. According to Reid's plurigus formula (\cite{R}), we have
$P_2(X)\ge 7$. This means ${\Cal O}(6)\hookrightarrow f_*\omega_{X'}^2$ and 
${\Cal E}\hookrightarrow f_*\omega_{X'}^4$. So $\fei{4}$ is generically finite, a contradiction. Finally, with regard to (A1), if $q(F)=0$, then $h^2(\Co{X'})=0$. So $\chi(\Co{X'})\le -1$. Thus $P_3(X)\ge 6$ according to Reid.  
We have ${\Cal O}(5)\hookrightarrow f_*\omega_{X'}^3$ and 
${\Cal E}\hookrightarrow f_*\omega_{X'}^5$. So $\fei{5}$ is generically finite, a contradiction. 
 
Step 2. $p_g(X)=1$ 

When $b=1$ or $b=0$ and $P_m(X)\ge 3$, we have $p_g(F)=1$ according to Proposition 3.1. This leads to (B0). {}From now on, we can suppose $b=0$, $P_m(X)=2$ and $p_g(F)\ge 2$.

We claim that $P_{m-2}(X)=1$. In fact, if $P_{m-2}(X)>1$, we must have $P_{m-2}(X)=2$. So the movable part of $|(m-2)K_{X'}|$ is a fiber $F$ of $f$. By Theorem 1.1, we have
$$|K_{X'}+\ulrcorner{\pi^*(K_X)}+F|\bigm|_F=|K_F+D|,$$ 
where $D:=\ulrcorner{\pi^*(K_X)}|_F$ is an effective divisor on $F$. So we see that $\dim\fei{m}(X)\ge 2$, a contradiction. 

If $P_{m-1}(X)=2$, then we can see from the above argument that $\dim\fei{m+1}(X)\ge 2$. This leads to (B3).

If $P_{m-1}(X)=1$ and we are not in (B1), then $P_{m+1}(X)=2$ by virtue of Proposition 3.1.  
We can easily see that $\dim\fei{m+2}(X)\ge 2$. This leads to (B2).

Step 3. $p_g(X)=0$ 

We can suppose $b=1$ and $q(F)\ge 2$ or $b=0$ and $q(F)\ge 3$. Otherwise we are in case (C0).

Suppose $b=1$. If $q(F)\ge 3$, we can see that $\dim\fei{m+1}(X)\ge 2$. In fact,
we can write 
$$m\pi^*(K_X)\sim_{\Bbb Q} aF+E_{\Bbb Q}^{(m)},$$
where $a=P_m(X)\ge 2$ and $E_{\Bbb Q}^{(m)}$ is an effective ${\Bbb Q}$-divisor. It is obvious that 
$$K_{X'}+\ulrcorner{m\pi^*(K_X)-F-\frac{1}{a}E_{\Bbb Q}^{(m)}}\le (m+1)K_{X'}.$$
Since $$m\pi^*(K_X)-F-\frac{1}{a}E_{\Bbb Q}^{(m)}\simnum m(1-\frac{1}{a})\pi^*(K_X)$$
is nef and big, we get by the vanishing theorem that
$$H^1(X', K_{X'}+\ulrcorner{m\pi^*(K_X)-F-\frac{1}{a}E_{\Bbb Q}^{(m)}})=0.$$
This means that
$$|K_{X'}+\ulrcorner{m\pi^*(K_X)-\frac{1}{a}E_{\Bbb Q}^{(m)}}|\bigm|_F=
\bigm| K_F+\ulrcorner{m\pi^*(K_X)-\frac{1}{a}E_{\Bbb Q}^{(m)}}|_F\bigm|,$$
where 
$$\ulrcorner{m\pi^*(K_X)-\frac{1}{a}E_{\Bbb Q}^{(m)}}|_F=\ulrcorner{(1-\frac{1}{a})E_{\Bbb Q}^{(m)}}|_F$$
is an effective divisor. According to \cite{X2}, $|K_F|$ gives a generically finite map. So we see that $\dim\fei{m+1}(F)=2$ and thus $\dim\fei{m+1}(X)\ge 2$. This leads to (C1).
  If $q(F)=2$ and $P_m(X)\ge 3$, we can still see that $\dim\fei{m+1}(X)\ge 2$. In this situation, $a\ge 3$. Let $F_1$ and $F_2$ be two distinct general fibers of $f$. Then we see that
$$ m\pi^*(K_X)-F_1-F_2-\frac{2}{a}E_{\Bbb Q}^{(m)}\simnum m(1-\frac{2}{a})\pi^*(K_X)$$
is nef and big. So we have the following surjective map 
$$\align
&H^0(X', K_{X'}+\ulrcorner{m\pi^*(K_X)-\frac{2}{a}E_{\Bbb Q}^{(m)}})\lrw  \\
&H^0(F_1, K_{F_1}+D_1)\oplus H^0(F_2, K_{F_2}+D_2)\lrw 0,
\endalign$$
where 
$$D_i=\ulrcorner{m\pi^*(K_X)-\frac{2}{a}E_{\Bbb Q}^{(m)}}|_{F_i}=
\ulrcorner{m(1-\frac{2}{a})E_{\Bbb Q}^{(m)}}|_{F_i}$$
is effective for all $i$.
This means that 
$$|K_{X'}+\ulrcorner{m\pi^*(K_X)-\frac{2}{a}E_{\Bbb Q}^{(m)}}|$$ 
can distinguish two different fibers of $f$ and 
$\dim\fei{m+1}(F_i)\ge 1$. We again see that $\dim\fei{m+1}(X)\ge 2$. This  leads to (C1).
If $q(F)=2$ and $P_m(X)=2$, we can use a parallel argument to that in the proof of the case $b=1$ of Step 1 to see that $\dim\fei{m+2}(X)\ge 2$. This corresponds to (C2).

Suppose $b=0$ and $q(F)\ge 3$. If $P_m(X)\ge 3$, we can use the same argument as in the case $b=1$ of Step 3 to see that $\dim\fei{m+1}(X)\ge 2$. This leads to (C1). What remains to be studied is the case $P_m(X)=2$. This is the most frustrating case. Anyway, it is easy to see that $\dim\fei{2m+1}(X)\ge 2$ in this case. Actually, one only has to consider the system 
$$|K_{X'}+\ulrcorner{m\pi^*(K_X)}+F|.$$
We can see that 
$$|K_{X'}+\ulrcorner{m\pi^*(K_X)}+F|\bigm|_F=\bigm|K_F+
\ulrcorner{m\pi^*(K_X)}|_F\bigm|,$$
where $\ulrcorner{m\pi^*(K_X)}|_F$ is effective.
This means 
$$\dim\fei{2m+1}(X)\ge\dim\fei{2m+1}(F)=2.$$
Now if $\dim\fei{2m}(X)\ge 2$, we are in (C3). If $\dim\fei{2m}(X)=1$, we have the following claim which shows that we are in either (C1) or (C4). We note that $P_{2m}(X)\ge 3$.
\medskip

\noindent {\bf Claim.} If $b=0$, $q(F)\ge 3$, $P_{m}(X)=2$, $P_{2m}(X)\ge 4$ and  
$\dim\fei{2m}(X)=1$. Then $\dim\fei{m+1}(X)\ge 2$. This leads to (C1).
\smallskip
 
Since $\dim\fei{2m}(X)=1$, we can see that both $\fei{2m}$ and $\fei{m}$ derive the same fibration $f:X'\lrw {\Bbb P}^1$. We can write
$$m\pi^*(K_X)\sim_{\Bbb Q} F+(D_v+D_h),$$
where $D_v$ is a vertical ${\Bbb Q}$-divisor with respect to the fibration $f$ and $D_h$ is the horizontal part. The supports of $D_v$ and $D_h$ are contained in the fixed part of $|mK_{X'}|$. $D_v$ and $D_h$ are both effective ${\Bbb Q}$-divisors. Similarly, we can write
$$2m\pi^*(K_X)\sim_{\Bbb Q} \sum_{i=1}^{a_2}F_i+(D_v'+D_h'),$$
where $a_2\ge P_{2m}(X)-1\ge 3$, $D_v'$ is a vertical effective ${\Bbb Q}$-divisor and $D_h'$ is a horizontal effective ${\Bbb Q}$-divisor. Since the support of 
$D_h'$ is contained in the fixed part of $|2mK_{X'}|$, we can see that 
$D_h'=2D_h$. Now we have
$$2m\pi^*(K_X)\simnum a_2F+D_v'+2D_h,$$
$$m\pi^*(K_X)\simnum \frac{a_2}{2}F+\frac{1}{2}D_v'+D_h.$$
So
$$m\pi^*(K_X)-F-\frac{1}{a_2}D_v'-\frac{2}{a_2}D_h\simnum m(1-\frac{2}{a_2})\pi^*(K_X)$$
is nef and big.  This means, according to the vanishing theorem, that
$$H^1(X', K_{X'}+\ulrcorner{m\pi^*(K_X)-F-\frac{1}{a_2}D_v'-\frac{2}{a_2}D_h})=0.$$
Denote $M:=\ulrcorner{m\pi^*(K_X)-\frac{1}{a_2}D_v'-\frac{2}{a_2}D_h}.$
Then $K_{X'}+M\le (m+1)K_{X'}$. Then
$$|K_{X'}+M|\bigm|_F=\bigm|K_F+M|_F\bigm|,$$
where $M|_F=\ulrcorner{(1-\frac{2}{a_2})D_h}|_F$ is an effective divisor on $F$.
So 
$$\dim\fei{m+1}(X)\ge\dim\Fi{K_{X'}+M}(F)=2.$$
The proof is complete.
\qed\enddemo

\proclaim{Theorem 3.4} Let $f: X\lrw C$ be a derived family of surfaces from the bicanonical pencil $|2K_X|$ of a smooth projective 3-fold $X$ of general type. 
Let $F$ be a general fiber of $f$ and denote $b:=g(C)$.
Then one of the following occurs.

(A0)' $p_g(X)>0$, $q(X)\le b+1$ and $p_g(F)=1$.

(A1)' $1\le p_g(X)\le 2$, $2\le P_2(X)\le 3$, $\dim\fei{3}(X)\ge 2$ and $p_g(F)\ge 2$.

(A2)' $p_g(X)=1$, $P_2(X)=P_3(X)=2$, $\dim\fei{4}(X)\ge 2$, $b=0$ and $p_g(F)\ge 2$.

(B0)' $p_g(X)=0$, $q(X)\le 2$ and  $q(F)\le 2$.

(B1)' $p_g(X)=0$, $\dim\fei{3}(X)\ge 2$ and $q(F)\ge 2$.

(B2)' $p_g(X)=0$, $P_2(X)=2$, $\dim\fei{4}(X)\ge 2$, $b=1$ and $q(F)=2$.

(B3)' $p_g(X)=0$, $P_2(X)=2$, $\dim\fei{4}(X)\ge 2$, $b=0$ and $q(F)\ge 3$.

(B4)' $p_g(X)=0$, $P_2(X)=2$, $P_4(X)=3$, $\dim\fei{5}(X)\ge 2$, $b=0$ and $q(F)\ge 3$.
\endproclaim
\demo{Proof} The proof is parallel to that of Theorem 3.3 except that we have more cases here. In order to avoid unnecessary redundancy, we only give the proof where it is different from the respective part in the proof of Theorem 3.3.

Step 1. $p_g(X)\ge 2$.

In this case, we always have $P_2(X)\ge 3$. When $b=0$ and $P_2(X)\ge 4$ or $b=1$, we see from Propositioin 3.1 that $p_g(F)=1$. This leads to (A0)'.

So we only have to consider the case with $b=0$, $P_2(X)=3$ and $p_g(F)\ge 2$.
In this situation, we see that the movable part of $|2K_{X'}|$ contains exactly 2 fibers of $f$. Since $P_3(X)\ge 4$, Proposition 3.1 gives  $\dim\fei{3}(X)\ge 2$. This corresponds to (A1)'.

Step 2. $p_g(X)=1$.

Excluding the situation (A0)' while observing Proposition 3.1, we only have to consider the case with $p_g(F)\ge 2$ and with the following extra properties:
$$b=0, \ P_2(X)\le 3\ \text{or}\ b=1,\ P_2(X)=2.$$

When $b=0$ and $P_2(X)=3$ or $b=1$ and $P_2(X)=2$, we know from Proposition 3.1 that $\dim\fei{3}(X)\ge 2$. This leads to (A1)'.

When $b=0$, $P_2(X)=2$, $p_g(F)\ge 2$ and $P_3(X)\ge 3$, we see from Proposition 3.1 that $\dim\fei{3}(X)\ge 2$. This also corresponds to (A1)'. Otherwise we always have $\dim\fei{4}(X)\ge 2$ because $P_4(X)\ge 3$. This is just (A2)'. 

Step 3. $p_g(X)=0$.

The argument in the proof of Theorem 3.3 is still effective in this case. We can see that (C0) through (C4) correspond to (B0)' through (B4)', respectively. We omit the proof.
\qed\enddemo

Now we can see that Theorem 1 (i), (iv) and (v), Theorem 2 are direct results from Theorem 3.3 and Theorem 3.4. In order to complete the proof of Theorem 1, we only have to show $q(X)\le 2$ whenever $m\ge 11$ or $m\ge 7$ and $p_g(X)>0$. 

\proclaim{Proposition 3.5} Let $X$ be a minimal projective 3-fold of general type with only ${\Bbb Q}$-factorial terminal singularities. Suppose $q(X)\ge 3$,
$P_{k_0}(X)>0$ and $P_{k_2}(X)\ge 2$. Then 
$\dim\fei{k_0+k_2+1}(X)\ge 2$. 
\endproclaim
\demo{Proof}
Choose a 1-dimensional subsystem $\Lambda\subset|k_2K_X|$ while taking a birational modification $\pi:X'\lrw X$ such that the pencil $\Lambda$ defines a morphism $g:X'\lrw {\Bbb P}^1$. We can even take further modification to $\pi$ so that $\pi^*(k_2K_X)$ has supports with only normal crossings. Taking the Stein factorization of $g$, then we get a derived fibration $p:X'\lrw C_1$. We note that this fibration is different from the one which was defined at the first page of this paper. Denote $b_1:=g(C_1)$. Let $M$ be the movable part of the pencil. We obviously have $M\le k_2K_{X'}$. We can write $M\simlin \sum_{i=1}^{a_1}F_i$, where $a_1\ge 1$ and $F_i$ is a fiber of $p$ for all $i$. 
We also note that $a_1=1$ if and only if $b_1=0$.
A general fiber $F$ is a smooth projective surface of general type. 

Suppose $b_1>0$. Then $|M|$ is base point free on $X$. Because $X$ has only isolated singularities, $F$ is smooth. We study the system $|tK_X+M|$ where $t\ge 2$.  We know that $M$ contains at least two components $F_1$ and $F_2$. By Theorem 1.2, we  see that 
$$H^0(X, tK_X+M)\lrw H^0(F_1, tK_{F_1})\oplus H^0(F_2, tK_{F_2})$$
is surjective.  This means that $\Fi{tK_X+M}$ can distinguish $F_1$ and $F_2$ and the restriction to $F_i$ is at least a bicanonical map. We know that $\dim\Fi{tK_{F_i}}(F_i)\ge 1$ for all $t\ge 2$. Noting that the image of $X$ through $\Fi{tK_X+M}$ is irreducible, we see that $\dim\Fi{tK_X+M}(X)\ge 2$.
So $\dim\fei{t+k_2}(X)\ge 2$. Thus $\dim\fei{k_0+k_2+1}(X)\ge 2$.

Suppose $b=0$. By Corollary 2.3, we have $q(F)\ge 3$. In this case, $M\simlin F$. We have
$$|K_{X'}+\ulrcorner{k_0\pi^*(K_X)}+F|\bigm|_F=\bigm|K_F+
\ulrcorner{k_0\pi^*(K_X)}|_F\bigm|,$$
where $\ulrcorner{k_0\pi^*(K_X)}|_F$ is effective. Thus $\dim\fei{k_0+k_2+1}(X)\ge\dim\fei{k_0+k_2+1}(F)=2$. 
\qed\enddemo

\proclaim{Lemma 3.6} Let $X$ be a smooth projective 3-fold of general type. If $q(X)\ge 3$, then either
$$P_2(X)>0\ \text{and}\ P_4(X)\ge 2$$
or 
$$P_k(X)\ge 2\ \text{for all}\ k\ge 5.$$
\endproclaim
\demo{Proof}
This is a byproduct from the proof of both Theorem 6.1, \cite{Ko1} and Proposition 4.3, \cite{Ko1}.  
\qed\enddemo

\proclaim{Proposition 3.7} Let $X$ be a minimal projective 3-fold of general type with only ${\Bbb Q}$-factorial terminal singularities. Suppose $q(X)\ge 3$, $P_2(X)>0$ and $P_4(X)\ge 2$. Then

(1) $\dim\fei{7}(X)\ge 2$. $\dim\fei{m}(X)\ge 2$ for all $m\ge 9$.

(2) If $\dim\fei{8}(X)=1$, then $p_g(X)=0$, $P_2(X)=1$, $P_4(X)=2$ and $P_8(X)=3$.
\endproclaim
\demo{Proof}  
(1) Let $k_0=2$ and $k_2=4$. Proposition 3.5 gives  $\dim\fei{7}(X)\ge 2.$
So $\dim\fei{2l+7}(X)\ge 2$ for all $l\in {\Bbb Z}^{+}$. 
Let $k_0=2$ and $k_2=7$. Applying Proposition 3.5 again, we get $\dim\fei{10}(X)\ge 2$. Thus $\dim_{2l+10}(X)\ge 2$ for all $l\in {\Bbb Z}^{+}$.

(2) We study $\fei{8}$. We have $P_8(X)\ge 3$. If $P_8(X)\ge 4$ and $\dim\fei{8}(X)=1$, we want to deduce a contradiction. We know that both $\fei{4}$ and $\fei{8}$ derive the same fibration $f:X'\lrw C$ which was described in the first page of this paper.  If $b>0$, it is easy to see that $\dim\fei{8}(X)\ge 2$ by a standard argument which has been used many times in this paper. So we can suppose $b=0$. Since $q(X)\ge 3$, we have $q(F)\ge 3$.
Suppose $M_4$, $M_8$ are the movable parts of $|4K_{X'}|$,  $|8K_{X'}|$ respectively. Then we have
$$4\pi^*(K_X)\sim_{\Bbb Q}M_4+E_4,$$
$$8\pi^*(K_X)\sim_{\Bbb Q}M_8+E_8,$$
where $E_4$ and $E_8$ are effective ${\Bbb Q}$-divisors. Let $E_v$, $E_v'$ be the vertical parts of $E_4$, $E_8$ and $E_h$, $E_h'$ be the horizontal parts   
of $E_4$, $E_8$ respectively. Because the support of $E_h$ is contained in the fixed part of $|4K_{X'}|$ and the support of $E_h'$ is contained in the fixed part of $|8K_{X'}|$, we see that $E_h'=2E_h$.
Now we have 
$$8\pi^*(K_X)\simnum a_8F+E_v'+2E_h,$$
where $a_8\ge 3$. It follows that
$$4\pi^*(K_X)\simnum \frac{a_8}{2}F+\frac{1}{2}E_v'+E_h.$$
Thus 
$$4\pi^*(K_X)-F-\frac{1}{a_8}E_v'-\frac{2}{a_8}E_h\simnum
4(1-\frac{2}{a_8})\pi^*(K_X)$$
is  a nef and big ${\Bbb Q}$-divisor. Denote 
$$G:=\ulrcorner{4\pi^*(K_X)-\frac{1}{a_8}E_v'-\frac{2}{a_8}E_h}.$$
Then we have 
$H^1(X', K_{X'}+G-F)=0$. So we see that
$$|K_{X'}+G|\bigm|_F=|K_F+G|_F|,$$
where 
$$G|_F=\ulrcorner{4\pi^*(K_X)-\frac{1}{a_8}E_v'-\frac{2}{a_8}E_h}|_F
=\ulrcorner{(1-\frac{2}{a_8})E_h}|_F$$
is effective.
So $\dim\fei{5}(X)\ge \dim\Fi{K_F+G|_F}(F)=2$. In particular, $P_5(X)\ge 2$. Now let $k_0=2$ and $k_2=5$. Applying Proposition 3.5, we see that 
$\dim\fei{8}(X)\ge 2$. This contradicts to our assumption. Thus we have seen that $P_8(X)=3$ if $\dim\fei{8}(X)=1$. It follows immediately that $P_4(X)=2$ and $P_2(X)=1$. If $P_g(X)>0$, it is very easy to see from Proposition 3.5 that $\dim\fei{8}(X)\ge 2$. 
So we have completed the proof.
\qed\enddemo

\proclaim{Proposition 3.8} Let $X$ be a minimal projective 3-fold of general type with only ${\Bbb Q}$-factorial terminal singularities. Suppose $q(X)\ge 3$ and $P_k(X)\ge 2$ for all $k\ge 5$. Then

(1) $\dim\fei{m}(X)\ge 2$ for all $m\ge 11$.

(2) If $\dim\fei{10}(X)=1$, then $p_g(X)=P_2(X)=P_3(X)=P_4(X)=0$.

(3) If $\dim\fei{9}(X)=1$, then $p_g(X)=P_2(X)=P_3(X)=0$.

(4) If $\dim\fei{8}(X)=1$, then $p_g(X)=P_2(X)=0$.

(5) If $\dim\fei{7}(X)=1$, then $p_g(X)=0$.

\endproclaim
\demo{Proof} 
Let $k_0=5$ and $k_2=t\ge 5$. Proposition 3.5 gives
$\dim\fei{t+6}(X)\ge 2$ for all $t\ge 5$. This leads to (1).

If $p_g(X)>0$, let $k_0=1$ and $k_2=t\ge 5$. Proposition 3.5 gives 
$\dim\fei{t+2}(X)\ge 2$ for all $t\ge 5$. 

If $P_2(X)>0$, let $k_0=2$ and $k_2=t\ge 5$. Proposition 3.5 gives 
$\dim\fei{t+3}(X)\ge 2$ for all $t\ge 5$.

If $P_3(X)>0$, let $k_0=3$ and $k_2=t\ge 5$. Proposition 3.5 gives 
$\dim\fei{t+4}(X)\ge 2$ for all $t\ge 5$.

If $P_4(X)>0$, let $k_0=4$ and $k_2=t\ge 5$. Proposition 3.5 gives 
$\dim\fei{t+5}(X)\ge 2$ for all $t\ge 5$.

(2), (3), (4) and (5) follow immediately.
\qed\enddemo

\proclaim{Corollary 3.9}  Let $f: X\lrw C$ be a derived family of surfaces from the m-canonical pencil $|mK_X|$ of a smooth projective 3-fold $X$ of general type.  Then

(1) $q(X)\le 2$ when $m\ge 11$.

(2) $q(X)\le 2$ when $m\ge 7$ and $p_g(X)>0$.
\endproclaim
\demo{Proof} This is obvious from Lemma 3.6, Proposition 3.7 and Proposition 3.8.
\qed\enddemo

\head 4. Appendix to Koll\'ar's method \endhead

Given a smooth projective 3-fold $X$ of general type, it is uncertain whether
$\dim\fei{m+1}(X)\ge \dim\fei{m}(X)$ for all $m>0$. Even if $P_m(X)>0$, it is false that $P_{m+1}(X)>0$. This makes it difficult to study some stable property of $\fei{m}$. That is why Koll\'ar's bound was bigger. Hereby we would like to study in an alternative way. The bounds are better, however still unsatisfactory. 

\proclaim{Proposition 4.1} Let $X$ be a minimal projective 3-fold of general type with only ${\Bbb Q}$-factorial terminal singularities. Suppose $q(X)\ge 2$ and $P_{k_2}(X)\ge 2$. Then $\dim\fei{m}(X)\ge 2$ for all $m\ge 4k_2+2$.
\endproclaim
\demo{Proof}
Choose a 1-dimensional subsystem $\Lambda\subset|k_2K_X|$ while taking a birational modification $\pi:X'\lrw X$ such that the pencil $\Lambda$ defines a morphism $g:X'\lrw {\Bbb P}^1$. We can even take further modification to $\pi$ so that $\pi^*(k_2K_X)$ has supports with only normal crossings. Taking the Stein factorization of $g$, then we get a derived fibration $p:X'\lrw C_2$. 
Denote $b_2:=g(C_2)$. Let $M$ be the movable part of the pencil. We obviously have $M\le k_2K_{X'}$. We can write $M\simlin \sum_{i=1}^{a_2}F_i$, where $a_2\ge 1$ and $F_i$ is a fiber of $q$. A general fiber $F$ is a smooth projective surface of general type. 

If $b_2>0$, then we can see that $\dim\fei{k_2+t}(X)\ge 2$ for all $t\ge 2$ according to the parallel argument in the proof of Proposition 3.5.

If $b_2=0$. we study in an alternative way. We have $M\simlin F$. Because 
$q(X)\ge 2$, we have $p_g(F)\ge q(F)\ge 2$. According to Theorem 1.1, we have
$$|K_{X'}+\ulrcorner{k_2\pi^*(K_X)}+F|\bigm|_F=\bigm|K_F+
\ulrcorner{k_2\pi^*(K_X)}|_F\bigm|.$$
This means $\dim\fei{2k_2+1}(F)\ge 1$ because $\ulrcorner{k_2\pi^*(K_X)}|_F$ is effective.  Suppose $M_{2k_2+1}$ is the movable part of $|(2k_2+1)K_{X'}|$ and 
$M_{2k_2+1}'$ is the movable part of 
$$|K_{X'}+\ulrcorner{k_2\pi^*(K_X)}+F|.$$
Then $M_{2k_2+1}'\le M_{2k_2+1}$. Let $M_0$ be the movable part of $|K_F|$. Then $h^0(F, M_0)\ge 2$. Considering the following two maps
$$H^0(X', K_{X'}+\ulrcorner{k_2\pi^*(K_X)}+F)\overset{\alpha}\to\lrw H^0(F, K_F+\ulrcorner{k_2\pi^*(K_X)}|_F)\lrw 0$$
$$H^0(X', M_{2k_2+1}')\overset{\beta}\to\lrw H^0(F, M_{2k_2+1}'|_F),$$
we know that $\alpha$ is surjective and the images of $\alpha$ and $\beta$ have the same dimension. So
$$\align
&h^0(F, M_{2k_2+1}'|_F)\ge \dim_{\Bbb C}\text{im}(\beta)=\dim_{\Bbb C}\text{im}(\alpha)\\
&=h^0(F,K_F+\ulrcorner{k_2\pi^*(K_X)}|_F).
\endalign$$
Because 
$$M_{2k_2+1}'|_F\le K_F+\ulrcorner{k_2\pi^*(K_X)}|_F,$$
we see that 
$$M_0\le M_{2k_2+1}'|_F\le M_{2k_2+1}|_F.$$
For all $t\ge 0$ and two different fibers $F_1$, $F_2$, we consider the system
$$|K_{X'}+\ulrcorner{(t+2k_2+1)\pi^*(K_X)}+F_1+F_2|.$$
It is obvious that
$$|K_{X'}+\ulrcorner{(t+2k_2+1)\pi^*(K_X)}+F_1+F_2|\subset |(t+4k_2+2)K_{X'}|.$$
{}From Theorem 1.1, we have the exact sequence
$$\align
&H^0(X',K_{X'}+\ulrcorner{(t+2k_2+1)\pi^*(K_X)}+F_1+F_2)\\
&\lrw H^0(F_1, K_{F_1}+G_1)\oplus H^0(F_2, K_{F_2}+G_2)\lrw 0,
\endalign$$
where $G_i=(\ulrcorner{(t+2k_2+1)\pi^*(K_X)}+F_1+F_2)|_{F_i}$ for all $i$.
We can see that
$$\align
K_{F_i}+G_i&\ge K_{F_i}+\ulrcorner{t\pi^*(K_X)|_{F_i}}+M_{2k_2+1}|_{F_i}\\
&\ge K_{F_i}+\ulrcorner{t\pi^*(K_X)|_{F_i}}+M_0
\endalign$$
for all $i$. Furthermore, one can see that
$$\dim\Fi{K_{F_i}+\ulrcorner{t\pi^*(K_X)|_{F_i}}+M_0}(F_i)\ge 1.$$
When $t=0$, it is obvious. When $t>0$, one need to use the vanishing theorem to prove it. Noting that the image of $X'$ through $\fei{t+4k_2+2}$ is irreducible and that 
$$\dim\fei{t+4k_2+2}(F_i)\ge 1$$ 
for all $i$, we can see $\dim\fei{t+4k_2+2}(X)\ge 2$.
\qed\enddemo

\proclaim{Proposition 4.2} Let $X$ be a minimal projective 3-fold of general type with only ${\Bbb Q}$-factorial terminal singularities. Suppose $q(X)>0$ and $P_{k_2}(X)\ge 2$. Then $\dim\fei{m}(X)\ge 2$ for all $m\ge 7k_2+3$.
\endproclaim
\demo{Proof}
We keep the same set up as in the proof of Proposition 4.1. We only have to study the case when $b_2=0$.  We have the fibration $p:X'\lrw {\Bbb P}^1$.
We still denote by $F$ a general fiber of $p$. Since $q(X)>0$, we get $p_g(F)\ge q(F)\ge 1$. If $p_g(F)\ge 2$, we have seen from the proof of the last proposition that we can get better bounds. The most frustrating case is when $p_g(F)=q(F)=1$. Let $\sigma: F\lrw F_0$ be the contraction onto the minimal model. According to Theorem 3.1 in \cite{Ci}, we know that $|2K_{F_0}|$ is base point free when $p_g(F)>0$. So the movable part of $|2K_F|$ is just $\sigma^*(2K_{F_0})$. According to Koll\'ar's method, we see that 
$$|(5k_2+2)K_{X'}|\bigm|_F\supset |2K_F|.$$
So, if we denote by $M_{5k_2+2}$ the movable part of $|(5k_2+2)K_{X'}|$,
we should have 
$$M_{5k_2+2}|_F\ge \sigma^*(2K_{F_0}).$$
For all $t\ge 0$ and two different fibers $F_1$, $F_2$, we consider the system
$$|K_{X'}+\ulrcorner{(t+5k_2+2)\pi^*(K_X)}+F_1+F_2|.$$
It is obvious that
$$|K_{X'}+\ulrcorner{(t+5k_2+2)\pi^*(K_X)}+F_1+F_2|\subset |(t+7k_2+3)K_{X'}|.$$
{}From Theorem 1.1, we have the exact sequence
$$\align
&H^0(X',K_{X'}+\ulrcorner{(t+5k_2+2)\pi^*(K_X)}+F_1+F_2)\\
&\lrw H^0(F_1, K_{F_1}+G_1')\oplus H^0(F_2, K_{F_2}+G_2')\lrw 0,
\endalign$$
where $G_i'=(\ulrcorner{(t+5k_2+2)\pi^*(K_X)}+F_1+F_2)|_{F_i}$ for all $i$.
We can see that
$$\align
K_{F_i}+G_i'&\ge K_{F_i}+\ulrcorner{t\pi^*(K_X)|_{F_i}}+M_{5k_2+2}|_{F_i}\\
&\ge K_{F_i}+\ulrcorner{t\pi^*(K_X)|_{F_i}}+\sigma^*(2K_{F_0})
\endalign$$
for all $i$.  Furthermore, one can see that
$$\dim\Fi{K_{F_i}+\ulrcorner{t\pi^*(K_X)|_{F_i}}+\sigma^*(2K_{F_0})}(F_i)\ge 1.$$
When $t=0$, it is obvious. When $t>0$, one need to use the vanishing theorem to prove it. Noting that the image of $X'$ through $\fei{t+7k_2+3}$ is irreducible and that 
$$\dim\fei{t+7k_2+3}(F_i)\ge 1$$ 
for all $i$, we can see $\dim\fei{t+7k_2+3}(X)\ge 2$.
The proof is complete.
\qed\enddemo

\proclaim{Corollary 4.3}  Let $f: X\lrw C$ be a derived family of surfaces from the m-canonical pencil $|mK_X|$ of a smooth projective 3-fold $X$ of general type. Then

(1) $q(X)\le 1$ whenever $m\ge 82$.

(2) $q(X)=0$ whenever $m\ge 143$.
\endproclaim
\demo{Proof} According to \cite{F} and Remark 6.6 in \cite{Ko1}, we always have $P_{20}(X)\ge 2$ if $q(X)>0$. 
 Let $k_2=20$ while applying Proposition 4.1 and Proposition 4.2, we get what we want. 
\qed\enddemo

\head Acknowledgment\endhead
This note was written while I was visiting as a post-doc fellow at the 
Mathematisches Institut der
Universit$\ddot{\text{a}}$t G$\ddot{\text{o}}$ttingen, Germany. I would like
to thank Prof. F. Catanese for useful discussions and helps during my stay
at G$\ddot{\text{o}}$ttingen.

\head {\bf References} \endhead
\roster
\item"[Be]" A. Beauville, {\it L'application canonique pour les surfaces de type g\'en\'eral}, Invent. Math. {\bf 55}(1979), 121-140.
\item"[Bo]" E. Bombieri, {\it Canonical models of surfaces of general type},
Publications I.H.E.S. {\bf 42}(1973), 171-219.
\item"[BPV]" W. Barth, C. Peter, A. Van de Ven, Compact complex surface, 1984, Springer-Verlag.
\item"[Ca1]" F. Catanese, {\it Canonical rings and ``special'' surfaces of general type}, Proc. Symposia in pure Math. {\bf 46}(1987), 175-194.
\item"[Ca2]" ------, {\it Singular bidouble covers and the construction of interesting algebraic surfaces}, Contemp. Math. {\bf 241}(1999), 97-120. 
\item"[Ci]" C. Ciliberto, {\it The bicanonical map for surfaces of general type}, Proc. Symposia in Pure Math. {\bf 62}(1997), 57-83.
\item"[E-L]" L. Ein, R. Lazarsfeld, {\it Global generation of pluricanonical and adjoint linear systems on smooth projective threefolds}, J. Amer. Math. Soc. {\bf 6}(1993), 875-903.
\item"[F]" A. R. Fletcher, {\it Contributions to Riemann-Roch in projective 3-folds with only canonical singularities and application}, Proc. Symposia in Pure Math. {\bf 46}(1987), 221-232.
\item"[H-N]" G. Harder, M. S. Narasimhan, {\it On the cohomology groups of moduli spaces of vector bundles on curves}, Math. Ann. {\bf 212}(1975), 215-248.
\item"[Ha]" R. Hartshorne, Algebraic Geometry, GTM {\bf 52}(1977),
Springer-Verlag.
\item"[Ka1]" Y. Kawamata, {\it A generalization of Kodaira-Ramanujam's
vanishing theorem}, Math. Ann. {\bf 261}(1982), 43-46.
\item"[Ka2]" ------, {\it Kodaira dimension of algebraic fiber space over curves}, Invent. Math. {\bf 66}(1982), 57-71.
\item"[KMM]" Y. Kawamata, K. Matsuda, K. Matsuki, {\it Introduction to the
minimal
model problem}, Adv. Stud. Pure Math. {\bf 10}(1987), 283-360.
\item"[Ko1]" J. Koll\'ar, {\it Higher direct images of dualizing sheaves}, I,
Ann. of Math. {\bf 123}(1986), 11-42.
\item"[Ko2]" ------, {\it Higher direct images of dualizing sheaves}, II, Ann. of Math. {\bf 124} (1986), 171-202.
\item"[K-M]" J. Koll\'ar, S. Mori, Birational geometry of algebraic
varieties, 1998, Cambridge Univ. Press.
\item"[N]" N. Nakayama, {\it Hodge filtrations and the higher direct images of canonical sheaves}, Invent. Math. {\bf 85}(1986), 237-251.
\item"[O]" K. Ohno, {\it Some inequalities for minimal fibrations of surfaces of general type over curves}, J. Math. Soc. Japan {\bf 44}(1992), 643-666.
\item"[R]" M. Reid, {\it Young person's guide to canonical singularities}, Proc. Symposia in pure Math. {\bf 46}(1987), 345-414.
\item"[V]" E. Viehweg, {\it Vanishing theorems}, J. reine angew. Math. {\bf
335}(1982), 1-8.
\item"[X1]" G. Xiao, {\it Finitude de l'application bicanonique des surfaces de type g\'en\'eral}, Bull. Soc. Math. France {\bf 113}(1985), 23-51.
\item"[X2]" ------, {\it L'irr\'egularit\'e des surfaces de type g\'en\'eral dont le syst\`eme canonique est compos\'e d'un pinceau}, Comp. Math. {\bf 56}(1985), 251-257.
\endroster
\enddocument